 \documentclass[11pt]{article}
\usepackage{amssymb, amsthm, amsmath, amscd}
\newtheorem{thm}{Theorem}[section]
\newtheorem{lem}[thm]{Lemma}

\newtheorem{cor}[thm]{Corollary}

\theoremstyle{definition}\newtheorem{df}[thm]{Definition}
\theoremstyle{definition}\newtheorem{rem}[thm]{Remark}
\theoremstyle{definition}

\renewcommand{\phi}{\varphi}

\newcommand{\C}{\mathbb{C}}
\newcommand{\T}{\mathbb{T}}

\newcommand{\Aff}{\operatorname{Aff}}

\newcommand{\morp}{contractive completely positive linear map}

\newcommand{\hm}{homomorphism}

\newcommand{\ep}{\epsilon}

\newcommand{\andeqn}{\,\,\,{\rm and}\,\,\,}
\newcommand{\rforal}{\,\,\,{\rm for\,\,\,all}\,\,\,}
\newcommand{\CA}{$C^*$-algebra}
\newcommand{\SCA}{$C^*$-subalgebra}

\newcommand{\beq}{\begin{eqnarray}}
\newcommand{\eneq}{\end{eqnarray}}
\newcommand{\tforal}{\,\,\,\text{for\,\,\,all}\,\,\,}

\title{ A short return to simple AH-algebras with real rank zero}
\author{Huaxin Lin
 }
\date{}

\begin{document}

\maketitle

\begin{abstract}
Let $A$ be a unital simple AH-algebra with stable rank one and real rank zero such that $kx=0$ for all $x\in
{\rm ker}\rho_A,$ the subgroup of infinitesmal
elements in $K_0(A),$ and for the same integer $k\ge 1.$  We show that $A$ has tracial rank zero and
is isomorphic to a unital simple AH-algebra with no dimension growth.

\end{abstract}

\section{Introduction}
One of the most successful aspects of operator algebras is the classification of simple separable amenable \CA s, or otherwise known as the Elliott program. The program started with the classification of unital simple
$A\T$-algebras with real rank zero up to isomorphisms by their ordered $K$-theory (with the scale) by G. A. Elliott \cite{Ell2} which was preceded
by Elliott's classification of AF-algebras (\cite{Ell1}) some fifteen years earlier.  It was followed  immediately by   a number
of earlier results. Then, Elliott and Gong (\cite{EG}) made a classification of AH-algebras with slow dimension growth and real rank zero by their scaled ordered $K$-theory. While the earlier results
concentrated in the cases of real rank zero, attention later shifted to the case that \CA s are not assumed to have real rank zero.
One of the highlights of the program is the classification (up to isomorphisms) of unital simple AH-algebras with no dimension growth
(see \cite{EGL}). These are \CA s whose real rank may not be zero (in fact
these \CA s have real rank one).
However, the Elliott invariant this time involves not only the ordered $K$-groups but also
tracial information.
Just as one thought a complete classification
for simple unital AH-algebras was possible, J. Villadsen provided examples
of unital simple AH-algebras whose stable rank may be greater than one and examples of unital simple AH-algebras with stable rank one whose $K_0$-group
may have perforation (\cite{V1} and \cite{V2}). It should be noted that unital simple AH-algebras
with slow dimension growth have stable rank one and have weakly unperforated $K_0$-groups.  On the other hand, it was proved in \cite{Lnproc} that a unital simple AH-algebra with stable rank one, real rank zero and with
weakly unperforated $K_0$-groups has tracial rank zero and therefore it is isomorphic to a unital simple AH-algebra with no dimension growth and with real rank zero. At that point, one might think the next goal would be the classification of unital simple AH-algebras of stable rank one and with weakly unperforated $K_0$-groups by their orderd $K$-groups and tracial information.  However, A. Toms provided (\cite{T})  an example of unital simple AH-algebra with stable rank one and with weakly unperforated $K_0$-group which is not isomorphic to a unital simple AH-algebra  with slow dimension growth and with the same Elliott invariant.

   The real rank of these AH-algebras are greater than zero.  The purpose of this note is to show that, when a unital simple AH-algebra has real rank zero, then it is more likely classifiable by its Elliott invariant.
Toms's first example is an inductive limit of homogeneous \CA s whose spectra are contractive finite CW complexes. As we will show in this short note, this example cannot be made so that it has real rank zero.  More precisely, we observe that a unital simple AH-algebra with real rank zero
which is an inductive limit of  homogeneous \CA s whose spectra are contractive finite CW complexes
is in fact  a unital simple AH-algebra with tracial rank zero. In particular, it is classifable  and isomorphic to a unital simple AH-algebra with no dimension growth.
Using our earlier result in \cite{Lnproc}, we further show  that if $A$ is a unital simple AH-algebra with stable rank one and real rank zero so that the subgroup of infinitesmal elements of $K_0(A)$ is a finite group (or zero), then $A$ is isomorphic to a unital simple AH-algebra with no dimension growth and therefore it is classifiable by the Elliott invariant. We actually prove a slightly more general result.

{\bf Acknowledgements}:  This work was done in the summer 2011 when the author was in East China Normal University. It supported by East China Normal University and the Changjiang Lectureship there. This work was also supported by a grant from NSF.

\section{The cases without torsion}

Let $A$ be a unital stably finite \CA. Denote by $T(A)$ the tracial state space of $A.$
 If $\tau\in T(A),$ we will also use $\tau$ for the trace
 $\tau\otimes Tr$ on $M_m(A),$ where $Tr$ is the standard trace on $M_m,$ where $m\ge 1$ is an integer. Let $\rho_A: K_0(A)\to \Aff(T(A))$ be the positive \hm\, defined by
$\rho_A([p])(\tau)=\tau(p)$ for all projections in $M_m(A),$ $m=1,2,....$

\begin{lem}\label{Ltr}
Let $A=\lim_{n\to\infty} (A_n, \phi_n)$ be a unital  \CA,
where $A_n$ is a unital \SCA\,  and $\phi_n: A_n\to A_{n+1}$ is a unital \hm. Suppose that
$p, q\in A$ are  two
non-zero projections and there is an integer $K\ge 1$
such that
$$
K\tau(q)<\tau(p) \rforal \tau\in T(A).
$$

Then,
there exist an integer $n(1)\ge 1,$ two
projections $p', q'\in A_{n(1)}$ and an integer $m_0\ge n(1)$ such that
$
\phi_{n(1), \infty}(p')
$
is unitarily equivalent to $p,$ $\phi_{n(1), \infty}(q')$ is unitarily equivalent to $q$
and
$$
Kt(\phi_{n(1), m}(q'))<t(\phi_{n(1),m}(p'))
$$
for all tracial states $t$ of $A_{m}$ with $m\ge m_0.$

\end{lem}

\begin{proof}
Put $B_n=\phi_{n, \infty}(A_n),$ $n=1,2,....$
Without loss of generality, we may assume
that there are $p', q'\in A_{n(1)}$ for some $n(1)$ such that
$$
\phi_{n(1), \infty}(p')=p\andeqn \phi_{n(1), \infty}(q')=q.
$$
Suppose that, for some increasing subsequence $\{n(k)\},$
\beq\label{Ltr-1-}
Kt_k(\phi_{n(1), n(k)}(q'))\ge t_k(\phi_{(n(1), n(k)}(p'))
\eneq
for some tracial states $t_k$ of $A_{n(k)}.$ By a result of Choi and Effros
(\cite{CE}), there exists a \morp\, $L_k: B_{k}\to A_{k}$ such that $\phi_{k,\infty}\circ L_k={\rm id}_{B_{k}},$ $k=1,2,....$
Note that, if $m\ge n(1),$
\beq\label{Ltr-1}
\lim_{k\to\infty}\| \phi_{m, n(k)}[L_m(\phi_{n(1), \infty}(b))-\phi_{n(1), m}(b)]\|=0
\eneq
for all $b\in A_{n(1)}.$ In particular, it holds for $m=n(2).$

Define $s_k': B_{n(k)}\to \C$ by $s_k'(b)=t_k\circ L_{n(k)}(b)$ for all $b\in B_{n(k)},$ $k=1,2,....$  Then $s_k'$ is a state on $B_{n(k)}.$ Let $s_k
$ be a state of $A$ which extends $s_k'.$  Let $\tau$ be a weak limit
of $\{t_k\}.$ It follows from (\ref{Ltr-1}) that
$\tau$ is a tracial state on $B_{n(1)}.$ One would have, by (\ref{Ltr-1-}),
$$
K\tau(q)\ge \tau(p).
$$
A contradiction.

\end{proof}

\begin{df}
{\rm
Recall that a unital \CA\, $A$ is said to be AH-algebra if
$A=\lim_{n\to\infty} (A_n, \phi_n),$ where $A_n=P_nM_{r(n)}(C(X_n))P_n,$
$X_n$ is a finite CW complex, $r(n)\ge 1$ is an integer and
$P_n\in M_{r(n)}(C(X_n))$ is a projection. We also assume that
$\phi_n: A_n\to A_{n+1}$ is unital. In what follows, we call $X_n$ the spectrum of $A_n.$

}
\end{df}

\begin{thm}\label{TCexmp}
Let $A=\lim_{n\to\infty}(A_n, \phi_n)$ be a unital simple AH-algebra
such that the spectrum of $A_n$ is a finite disjoint union of contractive spaces and $\phi_n: A_n\to A_{n+1}$ is a unital \hm. Suppose also
that $A$  has the real rank zero. Then
$A$ has tracial rank zero.
\end{thm}

\begin{proof}
Let $p, q\in A$ be two projections such that
$$
\tau(q)<\tau(p)
$$
for all $\tau\in T(A).$ We will show that
$q\lesssim p.$

We assume that $A_n=M_{r(n)}(C(X_n)),$ where $X_n$ is a disjoint union of contractive spaces.
We may also assume that there are two projections
$p', q'\in A_n$ for some $n\ge 1$ such that
$\phi_{n, \infty}(p')=p$ and $\phi_{n, \infty}(q')=q.$  Moreover,
by \ref{Ltr}, we may assume that
$$
t(q')<t(p')
$$
for all tracial states $t$ of $A_n.$
Since $X_n$ is a disjoint union of contractive spaces, all projections are trivial. Therefore
in $A_n,$
$q'\lesssim p'.$  It follows that
$q\lesssim p.$
This implies that $K_0(A)$ is weakly unperforated.

We will show that projections in $A$ have cancellation property.
Suppose that $p, q\in A$ and $p$ and $q$ are equivalent. As above,
there are $p', q'\in A_n$ such that
$\phi_{n, \infty}(p')=p$ and $\phi_{n, \infty}(q')=q.$ Moreover, we may assume that $p'$ and $q'$ are equivalent.
So in $A_n,$ $1-p'$ and $1-q'$ have the same rank at each point.
Note that, since $X_n$ is a disjoint union of contractive spaces,  $1-p'$ and $1-q'$ are trivial projections. It follows that $1-p'$ and $1-q'$ are equivalent.
It follows that $1-p$ and $1-q$ are equivalent in $A.$

Since $A$ is also assumed to have real rank zero, by (the proof of) part (3) of III.2.4 of \cite{BH}, $A$ has stable rank one.
It follows from
\cite{Lnproc} that $A$ has tracial rank zero.

\end{proof}

\begin{cor}\label{zerorho}
Let $A=\lim_{n\to\infty}(A_n, \phi_n)$ be a unital simple AH-algebra
such that ${\rm ker}\rho_{A_n}=\{0\}.$ Suppose that $A$ has real rank zero and stable rank one. Then $A$ has tracial rank zero and is isomorphic to a unital simple AH-algebra with no dimension growth.

\end{cor}

\begin{proof}
The proof is the same as that of \ref{TCexmp}. Note that since
${\rm ker}\rho_{A_n}=\{0\},$ two projections $p, q$ in $M_m(A_n)$ with the same rank must give the same element in $K_0(A_n).$
\end{proof}

\begin{lem}\label{proj}
Let $A$ be a unital AH-algebra and $x\in \rho_A(K_0(A))$ such that $x>0.$
Then, there is a projection $e\in M_m(A)$ for some integer $m\ge 1$ such that $\rho_A([e])=x.$
\end{lem}

\begin{proof}
Note that $M_m(A)$ is also a unital AH-algebra for any integer $m\ge 1.$
We may assume that there are two projections $p, q\in A$ such
that
$$
\rho([p]-[q])=x.
$$
We write $A=\lim_{n\to\infty}(A_n, \phi_n),$ where
$A_n=P_nM_{r(n)}(C(X_n))P_n,$ $X_n$ is a finite CW complex, and $r(n)\ge 1$ is an integer and $P_n\in M_{r(n)}(C(X_n))$ is a projection.
We may assume that there is an integer $n\ge 1$ and projections
$p', q'\in A_n$ such that
$\phi_{n, \infty}(p')=p$ and $\phi_{n, \infty}(q')=q.$ By \ref{Ltr},
we may also assume that
\beq\label{proj-1}
t(q')<t(p')
\eneq
for all tracial states $t$ of $A_n.$
There is an integer $K(n)\ge 1$  and there are trivial projections
$p_0, q_0\in M_{K(n)}(A_n)$ such that
$p_0$ has the same rank at each point of $X$ as that of $p',$ and
$q_0$ has the same rank at each point of $X$ as that of $q'.$ Therefore
\beq\label{proj-2}
t(p_0)=t(p')\andeqn t(q_0)=t(q')
\eneq
for all tracial states $t$ of $A_n.$ It follows that
\beq\label{proj-3}
\tau(\phi_{n, \infty}(p_0))=\tau(p)\andeqn
\tau(\phi_{n, \infty}(q_0))=\tau(q)
\eneq
for all $\tau\in T(A).$  On the other hand,
by (\ref{proj-1}) and (\ref{proj-2}), since both $p_0$ and $q_0$ are trivial, there is a partial isometry $v_0\in A_n$ such that
$$
v_0^*v_0=q_0\andeqn v_0v_0^*\le p_0.
$$
Let $v=\phi_{n, \infty}(v_0)$ and let $e=p_0-vv^*.$ Then
$e\in M_m(A)$ is a non-zero projection and, by (\ref{proj-3}),
$$
\rho_A([e])=x.
$$

\end{proof}

From the above and combining the result in \cite{Lnproc}, one has the following corollary:

\begin{cor}\label{C1}
Let $A$ be a unital simple AH-algebra with real rank zero and stable rank one. Suppose that ${\rm ker}\rho_A=\{0\}.$ Then $A$ has tracial rank zero.
\end{cor}

Note that, \ref{C1} is not a generalization of \ref{zerorho}. In fact
that a simple AF-algebra may have non-zero infinitesmal elements in
its $K_0$-group.

We will prove a much more general result that will allow  non-zero ${\rm ker}\rho_A.$

\section{The case that $\rho_A(K_0(A))$ is torsion}

\begin{lem}\label{KL1}
Let $A$ be a unital simple AH-algebra with stable rank one and with real rank zero. Let $K\ge 1$ be  an integer. Suppose that for any $x\in {\rm ker}\rho_A,$ $Kx=0.$
Suppose that $p, q\in M_m(A)$ are two projections
such that
$(K+2)\tau(q)<\tau(p)$
for all $\tau\in T(A).$
Then
$$
[q]\le [p].
$$

\end{lem}

\begin{proof}
Since $A$ has real rank zero, by a result of S. Zhang (\cite{Zh}, see also 9.4 of \cite{Lnmem}),
there exist mutually orthogonal projections
$e_1,e_2,...,e_{K+1}\in pAp$ such that
$e_i$ is equivalent to $e_1$ for $i=1,2,...,K,$  $e_{K+1}\lesssim e_1$
and $p=\sum_{i=1}^{K+1}e_i.$
Write $A=\lim_{n\to\infty}(A_n, \phi_n),$ where $A_n=P_nM_{r(n)}(C(X))P_n,$
$X$ is a finite CW complex, $r(n)\ge 1$ is an integer  and $P_n
\in M_{r(n)}(C(X))$ is a projection. Without loss of generality, we may assume that there are projections
$$
p', q', e_1',e_2',...,e_{K+1}'\in M_m(A_n)
$$
such that $\phi_{n, \infty}(p')=p,$
$\phi_{n, \infty}(q')=q,$ $\phi_{n, \infty}(e_i')=e_i,$
$i=1,2,...,K+1.$ We may also assume that
$e_i'\le p',$ $i=1,2,...,K+1,$ $e_1',e_2',...,e_{K+1}'$ are mutually orthogonal, $e_i'$ is equivalent to $e_1'$ in $M_m(A_n),$ $i=1,2,...,K,$
$e_{K+1}'\lesssim e_1'$ in $M_m(A_n)$
and $\sum_{i=1}^{K+1}e_i'=p'.$ Moreover, by \ref{Ltr}, we may also assume that
$$
(K+2)t(q')<t(p')\le (K+1)t(e_1')
$$
for all tracial states $t$ of $A_n.$
Note also, we have
\beq\label{KL-1}
(K+1)t(q')<Kt(e_1')
\eneq
for all tracial states $t$ of $A_n.$

There are trivial projections
$q'', e_1''\in M_R(A_n)$ for some $R\ge 1$ such that
$$
t(q'')=t(q')\andeqn t(e_1'')=t(e_1')
$$
for all tracial states $t\in T(A_n).$
Let ${\bar q}, {\bar e}\in M_R(A)$ such that
$\phi_{n, \infty}(q'')={\bar q}$ and $\phi_{n, \infty}(e_1'')={\bar e}.$
It follows that
$$
\tau({\bar q})=\tau(q)\andeqn \tau({\bar e})=\tau(e_1)
$$
for all $\tau\in T(A).$
Therefore
$$
[{\bar q}]-[q]\andeqn [{\bar e}]-[e_1]\,\,\,{\rm are\,\,\, in}\,\,\, {\rm ker}\,\rho_A.
$$
It follows from the assumption that
\beq\label{KL-2}
K[q]=K[{\bar q}]\andeqn K[{\bar e}]=K[e_1].
\eneq
On the other hand, by (\ref{KL-1}), since both $q''$ and $e_1''$ are trivial,
\beq\label{KL-3}
q''\lesssim e_1''.
\eneq
By (\ref{KL-2}) and (\ref{KL-3}),
\beq\label{KL-4}
K[q]=K[{\bar q}]\le K[e_1]\le p.
\eneq
Since $A$ also has stable rank one,
\beq\label{KL-5}
q\lesssim e_1+e_2+\cdots e_K\le p.
\eneq

\end{proof}

The following theorem was proved in \cite{Lnproc}.

\begin{thm}\label{Ttrtr0}
Let $A$ be a unital simple AH-algebra with stable rank one and real rank zero. Then, for any $\ep>0,$ $\sigma>0$ and any finite subset ${\cal F}\subset A,$
there exists a non-zero projection $p\in A$ and a finite dimensional \SCA\, $B\subset A$ with $1_B=p$ such that
\beq\label{Ttr0-1}
\|px-xp\|<\ep\tforal x\in {\cal F},\\
{\rm dist}(pxp, B)<\ep\tforal x\in {\cal F}\andeqn\\
\tau(1-p)<\sigma\tforal \tau\in T(A).
\eneq
\end{thm}

\begin{cor}\label{Ctr}
Let $A$ be a unital simple AH-algebra with stable rank one and real rank zero. Let $e_1, e_2\in A$ be two projections such that
\beq\label{Ctr0-1-1}
\tau(e_1)>\tau(e_2)\tforal \tau\in T(A).
\eneq
Then, for any $\ep>0,$ $\sigma>0$ and any finite subset ${\cal F}\subset A,$
there exists a non-zero projection $p\in A$ and a finite dimensional \SCA\, $B\subset A$ with $1_B=p$ such that
\beq\label{Ctr0-1}
\|px-xp\|<\ep\tforal x\in {\cal F},\\
{\rm dist}(pxp, B)<\ep\tforal x\in {\cal F}\andeqn\\
\tau(1-p)<\sigma\tforal \tau\in T(A).
\eneq
Moreover, there are projections $e_{j,0}\in (1-p)A(1-p)$ and
$e_{j,1}\in B,$ $j=1,2,$ such that
\beq\label{Ctr-1}
\|e_{j,0}+e_{j,1}-e_j\|<\min\{1/2, \ep\},\,\,\,j=1,2,\\
t(e_{1,1})> t(e_{2,1})\tforal  tracial \,\,\,states\,\,
of\,\,B\,\,\, and\\
\tau(e_{1,0})>\tau(e_{2,0})\tforal\tau\in T(A).
\eneq

\end{cor}

\begin{proof}
Note that $A$ is separable. Let $\ep>0,$ $\sigma>0$ and let
${\cal F}\subset A$ be a finite subset. Let
$\{x_1,x_2,...,\}$ be a dense subset of $A.$
Let ${\cal F}_n={\cal F}\cup\{e_1, e_2\} \cup \{x_1, x_2,..., x_n\}$ and let $\ep_n=\ep/2^n$ and $\sigma_n=\sigma/2^n,$ $n=1,2,....$ By \ref{Ttrtr0},
there exists a sequence of projections $p_n\in A$ and finite dimensional
\SCA s $B_n\subset A$ with $1_{B_n}=p_n$ such that
\beq\label{Ctr-3}
\|p_nx-xp_n\|<\ep/2^n\tforal x\in {\cal F}_n,\\
{\rm dist}(p_nxp_n, B_n)<\ep/2^n\tforal x\in {\cal F}_n\andeqn
\tau(1-p_n)<\sigma/2^n
\eneq
$n=1,2,....$
For all sufficiently large $n,$ there are projections
$e^{(j,1,n)}\in B_n$ and $e^{(j,0,n)}\in (1-p_n)A(1-p_n)$ such that
\beq\label{Ctr-4}
\|e^{(j,1,n)}+e^{(j,0,n)}-e_j\|<\ep/2^{n-2},
\eneq
$j=1,2$ and $n=1,2,....$
We also have
\beq\label{Ctr-4+}
\lim_{n\to\infty}\|p_ne_jp_n-e^{(j,1,n)}\|=0,\,\,\,j=1,2.
\eneq
Suppose that, for a subsequence $\{n_k\},$
\beq\label{Ctr-5}
t_k(e^{(2,1,n_k)})\le t_k(e^{(1,1,n_k)})
\eneq
for some tracial state $t_k$ of $B_{n_k}.$
Let $s_k': pAp\to \C$ be a state which extends $t_k.$
Define $s_k: A\to \C$ by $s_k(a)=s_k'(pap)$ for all $a\in A.$ Then
$s_k$ is a state on $A.$ Let $t$ be a weak limit of $\{s_k\}.$
One checks that $t$ is a tracial state on $A.$ Then
(\ref{Ctr-5}) and (\ref{Ctr-4+}) imply that
$$
t(e_1)\le t(e_2)
$$
which contradicts with (\ref{Ctr0-1-1}).

It follows that, for all sufficiently large $n,$
\beq\label{Ctr-6}
t(e^{(2,1,n)})>t(e^{(1,1,n)})\rforal t\in T(B_n).
\eneq
The lemma follows.

\end{proof}

\begin{thm}\label{MT}
Let $A$ be a unital simple AH-algebra with stable rank one and real rank zero. Suppose that there is an integer $K\ge 1$ such that, for any $x\in {\rm ker}\rho_A,$ $Kx=0.$ Then $A$ has tracial rank zero. Moreover, $A$ is isomorphic to a unital simple AH-algebra with slow dimension growth.
\end{thm}

\begin{proof}
We may assume that $A$ is infinite dimensional.
We will show that $K_0(A)$ is weakly unperforated. It then follows from
\cite{Lnproc} that $A$ has tracial rank zero.

It suffices to show the following:
If $p, q\in M_m(A)$ are two non-zero projections for some integer $m\ge 1$
and
$$
\tau(p)>\tau(q)\tforal \tau\in T(A),
$$
then $q\lesssim p.$

To prove this,
we note that  $M_m(A)$ is also a unital simple AH-algebra with stable rank one,
real rank zero  and $K_0(M_m(A))=K_0(A),$ so, to simplify the notation, without loss of generality,
we may assume that $p, q\in A.$

Let
$$
d_1=\inf\{\tau(p)-\tau(q): \tau\in T(A)\}.
$$
Since $A$ is simple, $d_1>0.$
Since  $A$ is an infinite dimensional simple \CA\, with real rank zero,
$pAp$ is also an infinite dimensional simple \CA\, with real rank zero.
It follows that there is a non-zero projection $e\in pAp$ such that
\beq\label{MT-0+1}
\tau(e)<d_1/2\tforal \tau\in T(A).
\eneq
Put
\beq\label{MT-0+d2}
d_2=\inf\{\tau(e): \tau\in T(A)\}.
\eneq
Note that $d_2>0.$
Put $p_0=p-e.$  Then
\beq\label{MT-0}
\tau(p_0)>\tau(q)\tforal \tau\in T(A).
\eneq
Let
$$
{\cal F}=\{p, q, e, p_0\}.
$$
It follows from \ref{Ttrtr0}  and \ref{Ctr} that, there exists a projection
$E\in A$ and a finite dimensional \SCA\, $B$ with $1_B=E$ such that
\beq\label{MT-2}
\|Ex-xE\|<{d_2\over{64K}}\tforal x\in {\cal F};\\\label{MT-2+1}
{\rm dist}(ExE, B)<{\min\{d_2, 1\}\over{64K}}\tforal x\in {\cal F}\andeqn\\\label{MT-2+2}
\tau(1-E)<{d_2\over{64K}}\tforal \tau\in T(A).
\eneq
Moreover, there are projections
$p_{0,1}, q_1\in B$ and $p_{0,0}, q_0\in (1-E)A(1-E)$ such that
\beq\label{MT-3}
\|p_{0,1}+p_{0,0}-p_0\|<{1\over{16K}},\\\label{MT-3+0}
\|q_1+q_0-q\|<{1\over{16K}}\andeqn\\\label{MT-3+1}
t(q_1)<t(p_{0,1})
\eneq
for all tracial state $t$ of $B.$
It follows that, in $B,$
\beq\label{MT-4}
q_1\lesssim p_{0,1}.
\eneq
We compute, by (\ref{MT-0+d2}) and (\ref{MT-2+2}) that
\beq\label{MT-5}
(K+2)\tau(q_0)<\tau(e)\tforal \tau\in T(A).
\eneq
It follows from \ref{KL1} that
\beq\label{MT-6}
q_0\lesssim e.
\eneq
Combining (\ref{MT-4}) and (\ref{MT-6}), we obtain that
$$
q_1+q_0\lesssim p_{0,1}+e.
$$
But, by (\ref{MT-3}) and (\ref{MT-3+0}),
$$
[q_1+q_0]=[q]\andeqn p_{0,1}\le p_0.
$$
Therefore
$$
q\lesssim p.
$$
\end{proof}

\begin{cor}\label{Mcor}
Let $A$ be a unital simple AH-algebra with stable rank one and real rank zero. Suppose that  ${\rm ker}\rho_A$ is finite. Then $A$ has tracial rank zero and
$K_0(A)$ is weakly unperforated.
\end{cor}

\section{Concluding remarks}

\begin{rem}
Theorem \ref{TCexmp} shows that that Toms' example (as in \cite{T})
could not occur under the assumption that $A$ has real rank zero.

\end{rem}

\begin{rem}
More can be said in Theorem \ref{TCexmp}. It is clear that that, in Theorem \ref{TCexmp}, it suffices to assume that every $A_n$ has the property that
all projections in the matrix algebras of $A_n$ are unitarily equivalent to those constant projections. So it allows $A_n$ to have  non-zero $K_1$-groups.

Let $A_n=P_nM_{r(n)}(C(X_n))P_n,$ where $X_n$ is a finite CW complex, $r(n)\ge 1$ is an integer and $P_n\in M_{r(n)}(C(X_n)).$ Suppose that $Y_n$ is another finite CW complex with covering dimension $hd(X_n)$ which has the same homotopy type of that of $X_n.$  Suppose
that
 $$
 \liminf_{n\to\infty} (\sup_{x\in X_n}{hd(X_n)\over{{\rm Rank} P_n(x)}})=0.
 $$
 Suppose also that $A=\lim_{n\to\infty}A_n$ is a unital simple \CA\, with real rank zero. Then, from the proof of \ref{TCexmp}, one can show that $A$ has tracial rank zero and stable rank one.  In other words, if $A$ is homotopically slow dimension growth and is of real rank zero and stable rank one, then $A$ is isomorphic to a unital simple AH-algebra with no dimension growth.
 One should note that, as in \cite{T}, without the assumption of real rank zero, the Cuntz semigroup of $A$ could be very different from those of unital simple AH-algebras with slow dimension growth.

\end{rem}

\begin{rem}
From the proof of \ref{MT}, one sees that Theorem \ref{MT} holds if
the assumption on ${\rm ker}\rho_A$ is replaced by the conclusion of \ref{KL1}, i.e., there is an integer $K\ge 1$ such that for any pair of
projections $p, q\in M_m(A)$ (for any integer $m\ge 1$), $K\tau(p)\le \tau(q)$ for all $\tau\in T(A)$ implies $p\lesssim q.$

\end{rem}

\noindent
hlin@uoregon.edu
\end{document}